\documentclass[11pt]{article}

\usepackage{amsmath}
\usepackage{amssymb}
\usepackage{amsfonts}

\def\U{{\cal{U}}}
\def\E{{\cal{E}}}

\def\a{{\mathbf{a}}}
\def\p{{\mathbf{p}}}
\RequirePackage{theorem}

\begin{document}

\title{
Words,  Dyck paths, Trees, and Bijections}
\author{Helmut Prodinger\\ 
\small{
The John Knopfmacher Centre for Applicable Analysis and Number
  Theory}\\
\small{  Department of Mathematics, University of the Witwatersrand}\\
  \small{P. O. Wits, 2050 Johannesburg, South Africa}\\
\small{\texttt{helmut@gauss.cam.wits.ac.za}}\\
\small{\texttt{http://www.wits.ac.za/helmut/index.htm}}
}

\date{}
\maketitle

\thispagestyle{empty}

\rightline{\sf{Dedicated to Gabriel Thierrin}}
\rightline{\sf{on the occasion of his eightieth birthday}}

\begin{abstract} 
In \cite{BaDeFePi96} the concept of nondecreasing Dyck paths was introduced.
We continue this research by looking at it from the point of view of
words, rational languages, 
planted plane trees, and
continued fractions.
 We construct a bijection with 
planted plane trees of height
$\le 4$ and compute various statistics on trees that are the equivalents of
nondecreasing Dyck paths.
\end{abstract} 

\section*{Personal reminiscences about Gabriel Thierrin}

Gabriel Thierrin invited me to London, Ontario, for six weeks
in February/March 1982. My memories about that trip are still
very much alive since it was my first crossing of the atlantic
ocean, and I was a very inexperienced traveller at that time, 
and everything was new to me. Snow storms, arctic temperatures,
clear blue skies, frozen sidewalks! I had the opportunity
to visit Waterloo, to give a talk, and to see the Niagara falls,
and to learn that Canada has best proportion of great rock
bands versus total population. 

I am only guessing, but I think that Gabriel Thierrin
was the referee of \cite{Prodinger80} and was attracted
by the combination of formal languages and other mathematical
concepts, perhaps not the standard ones in this context. 
I always liked the concepts of words, languages, grammars,
and automata, but I also wanted to see them in a wider context,
mostly in a combinatorial one. This is still true today, 
where I only occasionally bump into some (formal) languages. 

Gabriel Thierrin invited me to his house several times, 
and he and his wife were extremely friendly and helpful.
Once, he gave a party, and David Borwein also attended.
Later in life I met his sons Jonathan and Peter. 

We also wrote the paper \cite{PrTh83} together.

I remember much more, more than about any other trip I guess,
but perhaps I should rather stop here.

\smallskip

In the technical part of this paper, I want to demonstrate
a charming interplay of Dyck paths (related to Dyck words, of
course), certain rational languages and their associated
generating functions (being best described as continued
fractions), and some families of trees. The form of the generating
functions cries out for bijections, and they are described in the
sequel. Several characteristic parameters are also counted.

\section{Introduction}

In the  paper \cite{BaDeFePi96}, the 
Italian authors come up with the lovely
new concept of {\sl nondecreasing Dyck paths.}
Dyck words are geometric renderings of Dyck
paths where an open bracket is coded by an upward step,
and a closing bracket by a downward step.
The condition ``nondecreasing'' means
roughly that the sequence of the altitudes of the
{\sl valleys\/} must be {\sl nondecreasing.} We prefer to think about it in 
terms of {\sl planted plane trees};
 there is an obvious and well--known
bijection,  \cite{BrKnRi72, SeFl96}.

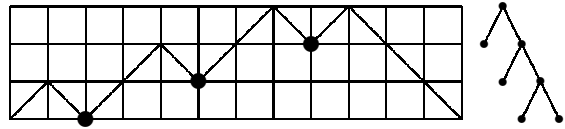
\begin{figure}[ht]  
\hskip2cm
\unitlength 0.50mm
\linethickness{0.4pt}
\begin{picture}(149.00,34.00)
\put(2.00,3.00){\line(1,0){120.00}}
\put(2.00,13.00){\line(1,0){120.00}}
\put(2.00,23.00){\line(1,0){120.00}}
\put(2.00,33.00){\line(1,0){120.00}}
\multiput(2.00,3.00)(0.12,0.12){84}{\line(0,1){0.12}}
\multiput(12.00,13.00)(0.12,-0.12){84}{\line(0,-1){0.12}}
\multiput(22.00,3.00)(0.12,0.12){167}{\line(0,1){0.12}}
\multiput(42.00,23.00)(0.12,-0.12){84}{\line(0,-1){0.12}}
\multiput(52.00,13.00)(0.12,0.12){167}{\line(0,1){0.12}}
\multiput(72.00,33.00)(0.12,-0.12){84}{\line(0,-1){0.12}}
\multiput(82.00,23.00)(0.12,0.12){84}{\line(0,1){0.12}}
\multiput(92.00,33.00)(0.12,-0.12){250}{\line(0,-1){0.12}}
\put(52.00,13.00){\circle*{4.00}}
\put(82.00,23.00){\circle*{4.00}}
\put(22.00,3.00){\circle*{4.00}}
\put(2.00,3.00){\line(0,1){30.00}}
\put(12.00,3.00){\line(0,1){30.00}}
\put(22.00,3.00){\line(0,1){30.00}}
\put(32.00,3.00){\line(0,1){30.00}}
\put(42.00,3.00){\line(0,1){30.00}}
\put(52.00,3.00){\line(0,1){30.00}}
\put(62.00,3.00){\line(0,1){30.00}}
\put(72.00,3.00){\line(0,1){30.00}}
\put(82.00,3.00){\line(0,1){30.00}}
\put(92.00,3.00){\line(0,1){30.00}}
\put(102.00,3.00){\line(0,1){30.00}}
\put(112.00,3.00){\line(0,1){30.00}}
\put(122.00,3.00){\line(0,1){30.00}}
\put(128.00,23.00){\circle*{2.00}}
\put(138.00,23.00){\circle*{2.00}}
\put(133.00,33.00){\circle*{2.00}}
\put(133.00,13.00){\circle*{2.00}}
\put(138.00,3.00){\circle*{2.00}}
\put(148.00,3.00){\circle*{2.00}}
\multiput(128.00,23.00)(0.12,0.24){42}{\line(0,1){0.24}}
\multiput(133.00,33.00)(0.12,-0.24){126}{\line(0,-1){0.24}}
\multiput(138.00,3.00)(0.12,0.24){42}{\line(0,1){0.24}}
\multiput(133.00,13.00)(0.12,0.24){42}{\line(0,1){0.24}}
\put(143.00,13.00){\circle*{2.00}}
\end{picture}
\caption{A nondecreasing Dyck path with valleys indicated
 and the corresponding
planted plane tree} 
\end{figure}
In honour of one of the authors, we decide to call the corresponding
trees {\sl Elena trees\/}, or simply 
{\sl Elenas\/}\footnote{In the literature, there are also 
{\sl Patricia\/} trees (tries).}.
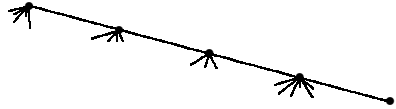
\begin{figure}[ht]  
\hskip2cm
\unitlength 0.50mm
\linethickness{0.4pt}
\begin{picture}(132.73,30.73)
\put(35.67,29.67){\circle*{2.11}}
\put(59.67,23.33){\circle*{2.11}}
\put(83.67,17.00){\circle*{2.11}}
\put(107.67,10.67){\circle*{2.11}}
\put(131.67,4.33){\circle*{2.11}}
\multiput(35.67,29.67)(0.45,-0.12){212}{\line(1,0){0.45}}
\put(131.67,4.33){\line(0,1){0.00}}
\multiput(30.33,28.33)(0.45,0.11){12}{\line(1,0){0.45}}
\multiput(35.67,29.67)(-0.20,-0.12){20}{\line(-1,0){0.20}}
\multiput(31.67,25.33)(0.12,0.14){31}{\line(0,1){0.14}}
\multiput(35.33,29.67)(-0.11,-0.45){6}{\line(0,-1){0.45}}
\multiput(35.67,29.67)(0.11,-2.00){3}{\line(0,-1){2.00}}
\multiput(52.33,21.00)(0.41,0.12){17}{\line(1,0){0.41}}
\multiput(59.33,23.00)(-0.12,-0.12){23}{\line(0,-1){0.12}}
\multiput(59.00,20.33)(0.11,1.00){3}{\line(0,1){1.00}}
\multiput(59.33,23.33)(0.11,-0.25){12}{\line(0,-1){0.25}}
\multiput(78.67,14.00)(0.19,0.12){26}{\line(1,0){0.19}}
\multiput(83.67,17.00)(-0.11,-0.31){12}{\line(0,-1){0.31}}
\multiput(83.67,16.67)(0.12,-0.22){17}{\line(0,-1){0.22}}
\multiput(101.00,8.67)(0.39,0.12){17}{\line(1,0){0.39}}
\multiput(107.67,10.67)(-0.15,-0.12){37}{\line(-1,0){0.15}}
\multiput(105.33,5.67)(0.12,0.25){20}{\line(0,1){0.25}}
\put(107.67,10.67){\line(0,-1){2.67}}
\multiput(111.00,5.33)(-0.12,0.19){26}{\line(0,1){0.19}}
\multiput(108.00,10.33)(0.14,-0.12){23}{\line(1,0){0.14}}
\end{picture}
\caption{A typical Elena; the short lines indicate paths of various
lengths} 
\end{figure}

In \cite{BaDeFePi96} the generating function of nondecreasing Dyck
paths of length $2n$ was already found to be
$\frac{z(1-z)}{1-3z+z^2}$. 
We find it practical also to include the {\sl empty path}, which gives us
\begin{equation*}
1+\frac{z(1-z)}{1-3z+z^2}=\frac{1-2z}{1-3z+z^2}\;.
\end{equation*}
Since the length $2n$ corresponds to an Elena of size ($=$number of nodes)
$n+1$, we find the generating function of Elenas as 
\begin{equation*}
E(z)=\sum_{n\ge0}[\text{\sf Number
of Elenas of size $n$}]\,z^n=\frac{z(1-2z)}{1-3z+z^2}\;.
\end{equation*}
Now here is an easy argument to see that directly.
We use the letter $\p$ to describe an arbitrary path of length $\ge1$
and the letter $\a$ which means 
`{\bf a}dvance to next node on the rightmost 
branch.' Then the set of Elenas $\cal{E}$
 is given by the symbolic equation (a rational language)
\begin{equation}\label{eq:elena}
\cal{E}=\big(\a\p^*\bigr)^*\a\;.
\end{equation}
Now mapping $\a\mapsto z$ and $\p\mapsto \frac{z}{1-z}$ we find
the generating function in the nice {\sl continued fraction\/} form
\begin{equation*}
E(z)={\cfrac{z}{1-\cfrac{z}{1-\cfrac{z}{1-z}}}}~.
\end{equation*} 

The continued fraction form suggests a relation to planted plane
trees of height $\le 4$; a bijection is constructed in the next section.

The following sections consider average values of several simple parameters
of Elenas. For simplicity, we give only first order asymptotics, but
explicit values (in terms of Fibonacci and Lucas numbers) and also
variances should not be too hard to obtain. 

Then we deal with the harder problem of the average height of (random)
Elenas of size $n$. 

We will use the number $\alpha=\frac{1+\sqrt5}{2}$
frequently in this paper, since it is prominent in the asymptotics
of Fibonacci numbers (and thus also Elenas).

\section{A bijection}

The continued fraction representation for $E(z)$ is 
well known in tree enumeration; it enumerates the set of planted
plane trees with height $\le 4$ (compare e.g. \cite{BrKnRi72,SeFl96}).

Now we will describe a bijection between Elenas and those trees.
We start from the representation $\big(\a\p^*\big)^*\a$
 and give an alternative
interpretation of the words in this set as height restricted trees.

First, a path with $n$  nodes (coded by $\p_n$) will be interpreted as
a root, followed by $n-1$ subtrees of size 1.

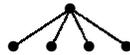
\begin{figure}[ht]  
\hskip2cm
\unitlength 0.50mm
\linethickness{0.4pt}
\begin{picture}(94.42,13.42)
\put(78.00,12.00){\circle*{2.83}}
\put(63.00,2.00){\circle*{2.83}}
\put(73.00,2.00){\circle*{2.83}}
\put(83.00,2.00){\circle*{2.83}}
\put(93.00,2.00){\circle*{2.83}}
\multiput(63.00,2.00)(0.18,0.12){84}{\line(1,0){0.18}}
\multiput(78.00,12.00)(-0.12,-0.24){42}{\line(0,-1){0.24}}
\multiput(83.00,2.00)(-0.12,0.24){42}{\line(0,1){0.24}}
\multiput(78.00,12.00)(0.18,-0.12){84}{\line(1,0){0.18}}
\end{picture}
\caption{Interpretation of a path with 5 nodes} 
\end{figure}

Then, a word $\a\p\dots \p$ will be interpreted as a root, followed
by subtrees given by the $\p$'s.  

\begin{figure}[ht]  
\hskip2cm
\unitlength 0.50mm
\linethickness{0.4pt}
\begin{picture}(135.00,24.42)
\put(25.00,3.00){\line(1,0){20.00}}
\put(45.00,3.00){\line(0,1){5.00}}
\put(45.00,8.00){\line(-1,0){20.00}}
\put(25.00,8.00){\line(0,-1){5.00}}
\put(55.00,3.00){\line(1,0){20.00}}
\put(85.00,3.00){\line(1,0){20.00}}
\put(115.00,3.00){\line(1,0){20.00}}
\put(75.00,3.00){\line(0,1){5.00}}
\put(105.00,3.00){\line(0,1){5.00}}
\put(135.00,3.00){\line(0,1){5.00}}
\put(75.00,8.00){\line(-1,0){20.00}}
\put(105.00,8.00){\line(-1,0){20.00}}
\put(135.00,8.00){\line(-1,0){20.00}}
\put(55.00,8.00){\line(0,-1){5.00}}
\put(85.00,8.00){\line(0,-1){5.00}}
\put(115.00,8.00){\line(0,-1){5.00}}
\put(80.00,23.00){\circle*{2.83}}
\multiput(35.00,8.00)(0.36,0.12){126}{\line(1,0){0.36}}
\multiput(80.00,23.00)(-0.12,-0.12){126}{\line(0,-1){0.12}}
\multiput(95.00,8.00)(-0.12,0.12){126}{\line(0,1){0.12}}
\multiput(80.00,23.00)(0.36,-0.12){126}{\line(1,0){0.36}}
\end{picture}
\caption{Interpretation of a $\a\p\p\p\p$; the boxes are the interpretations
of the respective paths } 
\end{figure}
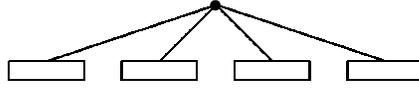

Finally, the last $\a$ will be the root, and the $\a\p\dots \p$'s become
subtrees of it.  
Figure~6 describes the process.

\begin{figure}[ht]  
\hskip2cm
\unitlength 0.50mm
\linethickness{0.4pt}
\begin{picture}(135.00,24.42)
\put(25.00,3.00){\line(1,0){20.00}}
\put(45.00,3.00){\line(0,1){5.00}}
\put(45.00,8.00){\line(-1,0){20.00}}
\put(25.00,8.00){\line(0,-1){5.00}}
\put(55.00,3.00){\line(1,0){20.00}}
\put(85.00,3.00){\line(1,0){20.00}}
\put(115.00,3.00){\line(1,0){20.00}}
\put(75.00,3.00){\line(0,1){5.00}}
\put(105.00,3.00){\line(0,1){5.00}}
\put(135.00,3.00){\line(0,1){5.00}}
\put(75.00,8.00){\line(-1,0){20.00}}
\put(105.00,8.00){\line(-1,0){20.00}}
\put(135.00,8.00){\line(-1,0){20.00}}
\put(55.00,8.00){\line(0,-1){5.00}}
\put(85.00,8.00){\line(0,-1){5.00}}
\put(115.00,8.00){\line(0,-1){5.00}}
\put(80.00,23.00){\circle*{2.83}}
\multiput(35.00,8.00)(0.36,0.12){126}{\line(1,0){0.36}}
\multiput(80.00,23.00)(-0.12,-0.12){126}{\line(0,-1){0.12}}
\multiput(95.00,8.00)(-0.12,0.12){126}{\line(0,1){0.12}}
\multiput(80.00,23.00)(0.36,-0.12){126}{\line(1,0){0.36}}
\put(26.00,4.00){\line(1,0){18.00}}
\put(44.00,4.00){\line(0,1){3.00}}
\put(44.00,7.00){\line(-1,0){18.00}}
\put(26.00,7.00){\line(0,-1){3.00}}
\put(56.00,4.00){\line(1,0){18.00}}
\put(86.00,4.00){\line(1,0){18.00}}
\put(116.00,4.00){\line(1,0){18.00}}
\put(74.00,4.00){\line(0,1){3.00}}
\put(104.00,4.00){\line(0,1){3.00}}
\put(134.00,4.00){\line(0,1){3.00}}
\put(74.00,7.00){\line(-1,0){18.00}}
\put(104.00,7.00){\line(-1,0){18.00}}
\put(134.00,7.00){\line(-1,0){18.00}}
\put(56.00,7.00){\line(0,-1){3.00}}
\put(86.00,7.00){\line(0,-1){3.00}}
\put(116.00,7.00){\line(0,-1){3.00}}
\end{picture}
\caption{Interpretation of a $(\a\p^*)(\a\p^*)(\a\p^*)(\a\p^*)\a$; 
the boxes are the interpretations
of the respective $(\a\p^*)$'s; the last $\a$ serves as the root} 
\end{figure}
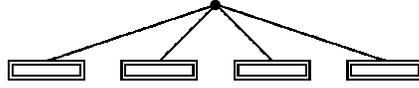

\begin{figure}[ht] \label{bijection}
\hskip2cm
\unitlength 0.50mm
\linethickness{0.4pt}
\begin{picture}(143.42,42.75)
\put(29.00,14.00){\circle*{2.83}}
\put(14.00,4.00){\circle*{2.83}}
\put(24.00,4.00){\circle*{2.83}}
\put(34.00,4.00){\circle*{2.83}}
\put(44.00,4.00){\circle*{2.83}}
\multiput(14.00,4.00)(0.18,0.12){84}{\line(1,0){0.18}}
\multiput(29.00,14.00)(-0.12,-0.24){42}{\line(0,-1){0.24}}
\multiput(34.00,4.00)(-0.12,0.24){42}{\line(0,1){0.24}}
\multiput(29.00,14.00)(0.18,-0.12){84}{\line(1,0){0.18}}
\put(59.00,14.00){\circle*{2.83}}
\put(54.00,4.00){\circle*{2.83}}
\put(64.00,4.00){\circle*{2.83}}
\multiput(59.00,14.00)(-0.12,-0.24){42}{\line(0,-1){0.24}}
\multiput(64.00,4.00)(-0.12,0.24){42}{\line(0,1){0.24}}
\put(74.00,14.00){\circle*{2.83}}
\put(83.00,4.00){\circle*{2.83}}
\put(93.00,4.00){\circle*{2.83}}
\put(103.00,4.00){\circle*{2.83}}
\put(119.00,14.00){\circle*{2.83}}
\put(114.00,4.00){\circle*{2.83}}
\put(124.00,4.00){\circle*{2.83}}
\multiput(119.00,14.00)(-0.12,-0.24){42}{\line(0,-1){0.24}}
\multiput(124.00,4.00)(-0.12,0.24){42}{\line(0,1){0.24}}
\put(131.00,14.00){\circle*{2.83}}
\put(142.00,14.00){\circle*{2.83}}
\put(92.67,14.33){\circle*{2.83}}
\multiput(83.00,4.00)(0.12,0.13){81}{\line(0,1){0.13}}
\multiput(92.67,14.33)(0.11,-3.44){3}{\line(0,-1){3.44}}
\multiput(92.67,14.33)(0.12,-0.12){87}{\line(1,0){0.12}}
\put(52.00,27.00){\circle*{2.83}}
\put(92.67,27.33){\circle*{2.83}}
\put(131.00,27.00){\circle*{2.83}}
\multiput(29.00,14.00)(0.21,0.12){109}{\line(1,0){0.21}}
\multiput(52.33,27.00)(0.12,-0.22){59}{\line(0,-1){0.22}}
\multiput(74.00,14.00)(-0.20,0.12){109}{\line(-1,0){0.20}}
\multiput(92.67,27.33)(0.11,-4.22){3}{\line(0,-1){4.22}}
\multiput(119.33,14.33)(0.12,0.13){100}{\line(0,1){0.13}}
\multiput(131.33,27.33)(-0.11,-4.44){3}{\line(0,-1){4.44}}
\multiput(142.33,14.00)(-0.12,0.14){92}{\line(0,1){0.14}}
\put(90.67,41.33){\circle*{2.83}}
\multiput(52.00,27.00)(0.33,0.12){120}{\line(1,0){0.33}}
\multiput(91.00,41.33)(0.12,-0.98){14}{\line(0,-1){0.98}}
\multiput(91.00,41.33)(0.34,-0.12){120}{\line(1,0){0.34}}
\put(110.67,27.33){\circle*{2.83}}
\multiput(110.67,27.33)(-0.17,0.12){114}{\line(-1,0){0.17}}
\end{picture}
\caption{Interpretation of  $(\a\p_5\p_3\p_1)(\a\p_4)(\a)
(\a\p_3\p_1\p_1)\a$; 
$\p_i$ stands for a path with $i$ nodes  } 
\end{figure}

Geometrically, one can imagine the process as follows. We consider
the rightmost branch of an Elena, take its right node as a root, move
the rest into horizontal position and rearrange the edges so that the nodes
are successors of the root. The attached paths are then rearranged as 
described.

\begin{figure}[ht]  
\hskip1cm
\unitlength 0.50mm
\linethickness{0.4pt}
\begin{picture}(186.00,66.00)
\put(86.00,6.00){\circle*{1.89}}
\put(86.00,11.00){\circle*{1.89}}
\put(86.00,16.00){\circle*{1.89}}
\put(86.00,21.00){\circle*{1.89}}
\put(107.00,10.00){\circle*{1.89}}
\put(112.00,15.00){\circle*{1.89}}
\put(117.00,20.00){\circle*{1.89}}
\put(122.00,15.00){\circle*{1.89}}
\put(181.00,10.00){\circle*{1.89}}
\put(176.00,15.00){\circle*{1.89}}
\put(171.00,20.00){\circle*{1.89}}
\put(166.00,15.00){\circle*{1.89}}
\put(56.00,12.00){\circle*{1.89}}
\put(61.00,17.00){\circle*{1.89}}
\put(66.00,12.00){\circle*{1.89}}
\put(137.00,10.00){\circle*{1.89}}
\put(142.00,15.00){\circle*{1.89}}
\put(147.00,10.00){\circle*{1.89}}
\put(142.00,20.00){\circle*{1.89}}
\put(61.00,12.00){\circle*{1.89}}
\multiput(137.00,9.67)(0.12,0.13){42}{\line(0,1){0.13}}
\put(142.00,15.00){\line(0,1){5.00}}
\multiput(142.00,15.00)(0.12,-0.12){42}{\line(0,-1){0.12}}
\multiput(55.67,12.00)(0.13,0.12){42}{\line(1,0){0.13}}
\put(61.00,17.00){\line(0,-1){5.00}}
\multiput(61.00,17.00)(0.12,-0.12){42}{\line(0,-1){0.12}}
\put(86.00,6.00){\line(0,1){15.00}}
\multiput(107.00,10.00)(0.12,0.12){84}{\line(0,1){0.12}}
\multiput(117.00,20.00)(0.12,-0.12){42}{\line(0,-1){0.12}}
\multiput(166.00,15.00)(0.12,0.12){42}{\line(0,1){0.12}}
\multiput(171.00,20.00)(0.12,-0.12){84}{\line(0,-1){0.12}}
\put(186.00,1.00){\line(0,1){25.00}}
\put(186.00,26.00){\line(0,1){25.00}}
\put(160.00,51.00){\line(0,-1){50.00}}
\put(128.00,51.00){\line(0,-1){50.00}}
\put(128.00,1.00){\line(0,1){0.00}}
\put(98.00,1.00){\line(0,1){50.00}}
\put(98.00,51.00){\line(0,1){0.00}}
\put(73.00,51.00){\line(0,-1){50.00}}
\put(48.00,51.00){\line(0,-1){50.00}}
\put(48.00,51.00){\line(0,1){15.00}}
\put(186.00,66.00){\line(0,-1){15.00}}
\put(160.00,66.00){\line(0,-1){15.00}}
\put(128.00,51.00){\line(0,1){15.00}}
\put(98.00,66.00){\line(0,-1){15.00}}
\put(73.00,51.00){\line(0,1){15.00}}
\put(61.00,31.00){\circle*{1.89}}
\put(61.00,36.00){\circle*{1.89}}
\put(61.00,41.00){\circle*{1.89}}
\put(61.00,46.00){\circle*{1.89}}
\put(61.00,31.00){\line(0,1){15.00}}
\put(77.00,35.00){\circle*{1.89}}
\put(82.00,40.00){\circle*{1.89}}
\put(87.00,45.00){\circle*{1.89}}
\put(92.00,40.00){\circle*{1.89}}
\multiput(77.00,35.00)(0.12,0.12){84}{\line(0,1){0.12}}
\multiput(87.00,45.00)(0.12,-0.12){42}{\line(0,-1){0.12}}
\put(123.00,35.00){\circle*{1.89}}
\put(118.00,40.00){\circle*{1.89}}
\put(113.00,45.00){\circle*{1.89}}
\put(108.00,40.00){\circle*{1.89}}
\multiput(108.00,40.00)(0.12,0.12){42}{\line(0,1){0.12}}
\multiput(113.00,45.00)(0.12,-0.12){84}{\line(0,-1){0.12}}
\put(138.00,38.00){\circle*{1.89}}
\put(143.00,43.00){\circle*{1.89}}
\put(148.00,38.00){\circle*{1.89}}
\put(143.00,38.00){\circle*{1.89}}
\multiput(137.67,38.00)(0.13,0.12){42}{\line(1,0){0.13}}
\put(143.00,43.00){\line(0,-1){5.00}}
\multiput(143.00,43.00)(0.12,-0.12){42}{\line(0,-1){0.12}}
\put(169.00,34.00){\circle*{1.89}}
\put(174.00,39.00){\circle*{1.89}}
\put(179.00,34.00){\circle*{1.89}}
\put(174.00,44.00){\circle*{1.89}}
\multiput(169.00,34.00)(0.12,0.12){42}{\line(0,1){0.12}}
\put(174.00,39.00){\line(0,1){5.00}}
\multiput(174.00,39.00)(0.12,-0.12){42}{\line(0,-1){0.12}}
\put(61.33,58.33){\makebox(0,0)[cc]{$\a\a\a\a$}}
\put(86.33,58.33){\makebox(0,0)[cc]{$\a\p_2\a$}}
\put(116.00,58.33){\makebox(0,0)[cc]{$\a\p_1\a\a$}}
\put(143.00,58.33){\makebox(0,0)[cc]{$\a\p_1\p_1\a$}}
\put(174.00,58.33){\makebox(0,0)[cc]{$\a\a\p_1\a$}}
\put(26.33,38.67){\makebox(0,0)[cc]{\small{Elena}}}
\put(26.33,18.33){\makebox(0,0)[cc]{\small{Height}}}
\put(26.33,11.33){\makebox(0,0)[cc]{\small{restricted}}}
\put(1.83,0.67){\line(0,1){51.25}}
\put(47.92,66.00){\line(1,0){132.92}}
\put(1.67,52.08){\line(1,0){179.17}}
\put(1.67,25.83){\line(1,0){179.17}}
\put(1.67,0.83){\line(1,0){179.17}}
\put(168.75,65.83){\line(1,0){17.08}}
\put(170.42,52.08){\line(1,0){15.42}}
\put(170.00,25.83){\line(1,0){15.83}}
\put(176.67,0.83){\line(1,0){9.17}}
\end{picture}
\caption{The bijection exemplified on trees with 4 nodes} 
\end{figure}

\section{Average degree of the root}

We use a second variable $u$ to label the degree of the root and
obtain easily
\begin{equation*}
T(z,u)=z+\cfrac{z}{1-\cfrac{uz}{1-z}}\frac{uz(1-2z)}{1-3z+z^2}\;.
\end{equation*}
To compute the average value, we have to
differentiate $T(z,u)$ with respect to $u$ and then to set $u=1$.
This yields
\begin{equation*}
\frac{\partial}{\partial u}T(z,u)\bigg|_{u=1}=
\frac{z^2(1-z)^2}{(1-2z)(1-3z+z^2)}\;.
\end{equation*} 
Around the (dominant) singularity $z=1/\alpha^2$ we have 
\begin{equation*}
\frac{z^2(1-z)^2}{(1-2z)(1-3z+z^2)}
\sim \frac{5-\sqrt5}{10}\frac{1}{1-z\alpha^2}\;,
\end{equation*}  
so that
\begin{equation*}
[z^n]\frac{z^2(1-z)^2}{(1-2z)(1-3z+z^2)}
\sim\frac{5-\sqrt5}{10}\alpha^{2n}~.
\end{equation*}  
Since
\begin{equation*}
[z^n]\frac{z(1-2z)}{1-3z+z^2}\sim \left(1-\frac{2}{\sqrt5}\right)
\alpha^{2n}
\end{equation*}
the average degree of the root is asymptotic to
\begin{equation*}
\frac{3+\sqrt{5}}{2}= 2.618033989
\;.
\end{equation*}

\section{Average number of leaves}

Replace $\a\mapsto z$  and $\p\mapsto\frac{zu}{1-z}$ 
in (\ref{eq:elena})
and multiply the whole thing by $u$ to get the
bivariate generating function
\begin{equation*}
\cfrac{zu}{1-\cfrac{z}{1-\cfrac{zu}{1-z}}}~.
\end{equation*}
Differentiate w.\;r.\;t. $u$, then set $u=1$ to get
$\dfrac{z(1-5z+8z^2-3z^3)}{(1-3z+z^2)^2}$.
Around the  singularity $z=1/\alpha^2$ we have 
\begin{equation*}
\frac{z(1-5z+8z^2-3z^3)}{(1-3z+z^2)^2}
\sim\frac{-2+\sqrt5}{5}\frac{1}{(1-z\alpha^2)^2}\;,
\end{equation*}  
so that
\begin{equation*}
[z^n]\frac{z(1-5z+8z^2-3z^3)}{(1-3z+z^2)^2}
\sim\frac{-2+\sqrt5}{5}n\alpha^{2n}\;.
\end{equation*}  
Since
\begin{equation*}
[z^n]\frac{z(1-2z)}{1-3z+z^2}\sim \left(1-\frac{2}{\sqrt5}\right)
\alpha^{2n}
\end{equation*}
the average number of leaves is asymptotic to
\begin{equation*}
\frac{n}{\sqrt5}= 0.4472135956\,n\;.
\end{equation*}

\section{Average number of paths}

Replace $\a\mapsto z$  and $\p\mapsto\frac{zu}{1-z}$ to get the
bivariate generating function
\begin{equation*}
\cfrac{z}{1-\cfrac{z}{1-\cfrac{zu}{1-z}}} ~.
\end{equation*} 
Differentiate w.\;r.\;t. $u$, then $u=1$ yields
$\dfrac{z^3(1-z)}{(1-3z+z^2)^2}$.
Hence \newline
$[z^n]\displaystyle{\frac{z^3(1-z)}{(1-3z+z^2)^2}
\sim\frac{-2+\sqrt5}{5}\,n\,\alpha^{2n}}
$.
Thus the average number of paths is
asymptotic to
\begin{equation*}
\frac{n}{\sqrt5}=0.4472135956\,n\;.
\end{equation*}

\section{Average number of nodes `$\a$' }

Replace $\a\mapsto zu$  and $\p\mapsto\frac{z}{1-z}$ to get the
bivariate generating function
\begin{equation*}
\cfrac{zu}{1-\cfrac{zu}{1-\cfrac{z}{1-z}}}~.
\end{equation*} 
Differentiate w.\;r.\;t. $u$, then $u=1$ yields
$\dfrac{z(1-2z)^2}{(1-3z+z^2)^2}$.
Hence

\begin{equation*}
[z^n]\frac{z(1-2z)^2}{(1-3z+z^2)^2}\sim\frac{7-3\sqrt5}{10}\,n\,\alpha^{2n}\;.
\end{equation*}
Thus the average number of $\a$'s  is
asymptotic to
\begin{equation*}
\frac{5-\sqrt5}{10}n= 0.2763932022\,n\;.
\end{equation*}
As a corollary, we get that the average number of nodes lying in {\sl paths\/} is
asymptotic to 
\begin{equation*}
n-\frac{5-\sqrt5}{10}n=\frac{5+\sqrt5}{10}n =
0.7236067978\,n\;.
\end{equation*}
And furthermore the average number of nodes in {\sl one\/} path is asymptotic to
\begin{equation*}
\frac{5+\sqrt5}{10}n\bigg/\frac{n}{\sqrt5}=\frac{1+\sqrt5}{2}= 1.618033989
\;.
\end{equation*}

\section{Number of descendants}

The number of descendants of a node is the size of the subtree with this node
as the root. 
The paper \cite{MaPaPr98} deals e.~g. extensively with this subject.
We want to know the average number of descendants. This is an
average over both, the Elenas, and the nodes in an Elena. Thus it is meaningful
to define for an Elena $t$ 
\begin{equation*}   
\psi(t):=\sum_{\text{$v$ a node of $t$  }}[\text{number of descendants of $v$}]
\end{equation*}
and
\begin{equation*}   
D(z,u):=\sum_{t\in \E}z^{|t|}u^{\psi(t)}\;;
\end{equation*}
then we find the desired average as
$\frac{1}{n}[z^n]\frac{\partial}{\partial u}D(z,u)\big|_{u=1}$.
Now we want to derive a functional equation for this function
$D(z,u)$. Of course we follow the general decomposition 
(\ref{eq:elena}).
The contribution of each path attached to the root is
\begin{equation*}   
Q(z,u)=\sum_{m\ge1}z^mu^{\binom{m+1}{2} }\;.
\end{equation*}
The contribution of the root is $zu^n$, which is handled by
first neglecting it and then substituting $zu$  for $z$. Altogether
we find
\begin{equation*}   
D(z,u)=zu+\frac{zu\,D(zu,u)}{1-Q(zu,u)}\;.
\end{equation*} 
Now let us differentiate this w.\;r.\;t. to $u$ and plug in $u=1$.
We can also use the special values
\begin{equation*}   
D(z,1)=E(z)\quad\text{and}\quad \frac{\partial }{\partial z }D(z,1)=
\frac{1-4z+5z^2}{(1-3z+z^2)^2}
\end{equation*}  
as well as
\begin{equation*}   
Q(z,1)= \frac{z}{1-z}\ \text{and}\ 
\frac{\partial }{\partial z }Q(z,1)=\frac{1}{(1-z)^2}
\ \text{and}\ 
\frac{\partial }{\partial u }Q(z,1)=\frac{z}{(1-z)^3}\;.
\end{equation*}
The resulting equation contains only one unknown function, $\frac{\partial}{\partial u}D(z,u)\big|_{u=1}$, and Maple solves it as
\begin{equation*}   
\frac{\partial}{\partial u}D(z,u)\bigg|_{u=1}=
\frac{z(1-7z+20z^2-26z^3+11z^4)}{(1-z)(1-3z+z^2)^3}
\sim\frac{7-3\sqrt5}{10}\frac{1}{(1-z\alpha^2)^3}
\;.
\end{equation*}
Hence
\begin{equation*}   
\frac{1}{n}[z^n]\frac{\partial}{\partial u}D(z,u)\bigg|_{u=1}\sim
\frac{7-3\sqrt5}{10}\,\frac{n}{2}\,\alpha^{2n}\;.
\end{equation*}
Dividing this quantity by the asymptotic
equivalent for the total number,
$\big(1-\tfrac{2}{\sqrt5}\big)\alpha^{2n}$, we get the average number
of descendants
as
\begin{equation*}   
\frac{5-\sqrt5}{20}\, n=  0.1381966011
\,n\;.
\end{equation*} 

\section{Number of ascendants}

The number of ascendants of a node is defined to be the number
of nodes on the path of the node to the root. It is also called
the {\sl depth}. And the sum over all depths (summed over all
nodes in the Elena) is called the {\sl path length}. 
It is very similar
to the {\sl area}, studied in the  paper \cite{BaDeFePi96}.

However, it is quite easy to see that the
average number of ascendants
equals the average number of descendants: Consider two nodes $i$ and $j$
such that $i$ lies on the path from the root to $j$. 
Then $i$ appears in the count of $j$ of the number of ascendants, and
$j$ appears in the count of $i$ of the number of descendants. Since these
quantities are summed over all nodes, we are done.
(This argument is general and not restricted to Elenas.)

\section{Average height of Elenas}

The recursion
$\E=a+(\a\p^*)\E$
translates into
$E=z+\frac{z(1-z)}{1-2z}E$
and also into the recursion for $E_h$, the generating functions of
Elenas of height $\le h$,
\begin{equation*}   
E_h=z+\frac{z(1-z)}{1-2z+z^h}E_{h-1}\;.
\end{equation*}
 Denoting the generating functions of Elenas of height $>h$ by
$U_h$, we find by taking differences
\begin{equation*}   
(1-2z+z^h)U_h=\frac{(1-z)z^{h+2}}{1-3z+z^2}+z(1-z)U_{h-1}\;.
\end{equation*} 
We find the average height as
\begin{equation*}   
\frac{[z^n]\sum_{h\ge0}U_h(z)}{[z^n]E(z)}\;.
\end{equation*}
Now define $\U(z,w):=\sum_{h\ge0}U_h(z)w^h$.
Summing up we get
\begin{multline*}   
(1-2z)\,\U(z,w)+\U(z,zw)\\
=\frac{2z(1-z)(1-2z)}{1-3z+z^2}
+\frac{z^3w(1-z)}{(1-3z+z^2)(1-wz)}+wz(1-z)\,\U(z,w)~.
\end{multline*}  
The instance $w=1$ is of special interest;
\begin{equation*}   
(1-3z+z^2)\,\U(z,1)+\U(z,z)=\frac{z(2-6z+5z^2)}{1-3z+z^2}\;.
\end{equation*}
 From this we see that $\U(z,1)$ has a double pole at the dominant singularity
$z=\lambda:=1/\alpha^2$. Since for all $h\geq0$ 
\begin{equation*}   
U_h(z)\sim\frac{\lambda(1-2\lambda)}{1-3z+z^2}~,
\end{equation*}
we infer that
$\U(z,z)\sim\dfrac{\lambda(1-2\lambda)}{1-\lambda}\dfrac1{1-3z+z^2}$.
Hence
\begin{equation*}   
\U(z,1)\sim\frac{47-21\sqrt5}{2}\frac{1}{(1-3z+z^2)^2}\sim\frac{7-3\sqrt5}{10}\frac{1}
{(1-z\alpha)^2}
\end{equation*}
and
\begin{equation*}   
[z^n]\,\U(z,1)\sim\frac{7-3\sqrt5}{10}\,n\,\alpha^{2n}\;.
\end{equation*}
Dividing this by $\big(1-\tfrac{2}{\sqrt5}\big)\alpha^{2n}$,
we find for the average height the asymptotic equivalent
\begin{equation*}   
\frac{5-\sqrt5}{10}n= 0.2763932022\,n\;.
\end{equation*}

\section{Conclusion}

For the reader's convenience we collect our findings in a small table.

\begin{table}
\renewcommand{\arraystretch}{4.0}
\begin{center}
\renewcommand{\arraystretch}{3.0}
\begin{tabular}[h]{| l | c |}
\hline
Degree of  root & $\dfrac{3+\sqrt5}{2}$ \\  \hline
Number of leaves & $\dfrac{n}{\sqrt5}$  \\ \hline
Number of paths &  $\dfrac{n}{\sqrt5}$\\ \hline
Number of nodes on rightmost branch & $\dfrac{5-\sqrt5}{10}{n}$  \\ \hline
Number of nodes in paths&  $\dfrac{5+\sqrt5}{10}{n}$\\ \hline
Number of nodes in one path&  $\dfrac{1+\sqrt5}2$\\ \hline
Number of ascendants&  $\dfrac{5-\sqrt5}{20}n$\\ \hline
Number of descendants&  $\dfrac{5-\sqrt5}{20}n$\\ \hline
Height&  $\dfrac{5-\sqrt5}{10}n$\\ \hline
\end{tabular}
\end{center}
\caption{Several averages on Elenas}
\end{table}

\bibliographystyle{plain}


\end{document}